\documentclass[12pt]{amsart}

\input psfig
\begin{document}

\title{ON A COMPUTER RECOGNITION OF 3-MANIFOLDS}

\author{Sergei V. Matveev}
\thanks{Research in MSRI was supported in part 
 by Russian Fond of Fundamental Investigations, grant N 96-01-847,
and by INTAS, project  94-921.
Research at MSRI is supported in part by NSF grant DMS-9022140}

\address{Sergei V. Matveev
Chelyabinsk State University
Chelyabinsk, 454136
Russia}

\email{matveev@cgu.chel.su}
   
\begin{abstract}
We describe theoretical backgrounds for a computer program that
recognizes all closed orientable 3-manifolds up to complexity 8. The program 
can treat also not necessarily closed 3-manifolds of  bigger complexities, but
here some unrecognizable (by the program) 3-manifolds may occur.
\end{abstract}
\maketitle

\section*{Introduction}

Let $M$ be an orientable 3-manifold such that $\partial M$ is either empty or
consists of  tori. Then, modulo W. Thurston geometrization conjecture [Scott
1983], $M$ can be decomposed in a unique way  into graph-manifolds 
and hyperbolic pieces. The classification of graph-manifolds is well-known
[Waldhausen 1967], and a list of cusped hyperbolic manifolds up to complexity 7 is
contained in [Hildebrand, Weeks 1989]. If we possess an information how the pieces are
glued together, we can get an explicit description  of $M$ as a sum of
geometric pieces. Usually such a presentation is sufficient for understanding
the intrinsic structure of $M$; it allows one to label $M$ with a {\it name }
that distinguishes it from all other manifolds.

We describe theoretical backgrounds and a general scheme of a computer algorithm
that realizes in part the procedure. Particularly, for all closed orientable 
manifolds up to complexity 8 (all of them are graph-manifolds, see
 [Matveev 1990])
 the algorithm gives an exact answer.

The paper is based on a talk at   MSRI workshop on computational and
algorithmic methods in three-dimensional topology (March 10-14, 1997). 
 The author wishes to thank MSRI for a friendly atmosphere
and good conditions of work.

\section{Special and almost special spines }  

Let $\Delta$ denote the underlying space of the one-dimensional skeleton of the 
3-simplex, i.e. the polyhedron homeomorphic to a circle with three radii.

{\bf Definition.} A compact polyhedron $P$ is called a {\it fake surface} if
the link of each of its points is homeomorphic to one of the following
one-dimensional polyhedra:

(1) a circle; (2) a circle with a diameter; (3) a circle with three radii (i.e.
$\Delta$).

Typical neighborhoods of points of a fake surface are shown in Fig. 1, where
these points are shown as fat dots. 

\begin{figure}
  \centerline{\psfig{file=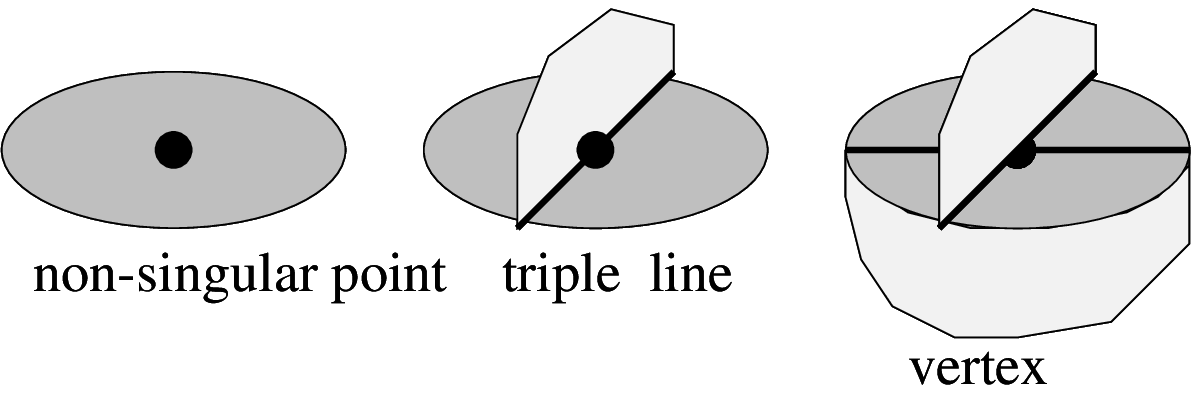,width=10cm}}  
  \caption{\label{details}}
  \end{figure}

The union of singular points (i.e. 
vertices and triple lines) of a fake surface $P$ is called its {\it singular 
graph} (and is denoted by $SP$). Connected components of $P\setminus SP$ are called 
{\it 2-components} of $P$. Each 2-component is a 2-manifold without boundary.

{\bf Definition.} A compact polyhedron $P$ is called  {\it almost special} if
it can be embedded in a fake surface.

There is a close relation between fake surfaces and almost special polyhedra.
For example, the wedge of any fake surface and any graph is an almost 
special polyhedron. The example is very typical, since any almost special
polyhedron can be collapsed onto a polyhedron of the form $F\cup G$, where
$F$ is a collection of disjoint fake surfaces, $G$ is a graph, and $F\cap G$ 
is a finite set of non-vertex points in $F$.

{\bf Definition.} A fake surface $P$ is said to be a {\it special polyhedron} if
it contains at least one vertex and if all its 2-components are 2-cells.

Note that there are finitely many special polyhedra with a given number of
vertices.

{\bf Definition.} A subpolyhedron $P\subset $ Int $M$ of a compact 3-manifold $M$ with
a non-empty boundary is said to be its {\it spine} if $M$ collapses to $ P$ or, equivalently, if 
$M\setminus P$ is homeomorphic to $\partial M\times (0,1].$

The spine is called almost special or special if it is a polyhedron
of the corresponding type.

We always assume that an almost special spine can not be collapsed onto a 
proper subpolyhedron.

By a spine of a closed manifold $M$ we  mean a spine of the punctured
$M$, i.e. $M\setminus $ Int $ D^3$. It is known that any 
compact 3-manifold possesses an
almost special (even a special) spine. Moreover, one can easily construct a 
special spine of $M$ starting from practically any presentation of $M$ 
[Matveev 1990]. Special spines possess an important property that favorably 
distinguishes them from fake surfaces and almost special spines: a 3-manifold
can be uniquely recovered from its special spine. Note that special spines
can be considered as combinatorial objects and admit presentations in
computer's memory as strings of integer numbers that show how 2-cells are
attached to singular graphs of spines. To present a manifold by its almost special
spine, an additional information is needed about the way how the spine should be 
thickened to the 3-manifold.

\section{Simplifying moves on spines}

In what follows we will consider compact orientable 3-manifolds whose boundaries
consist of spheres and tori.

 We introduce five types of moves on almost special spines. Each of
them does not increase the number of vertices of the spine and quite often
decreases it. We call them {\it simplifying moves}. Any spine admits only 
finitely many simplifying moves. The moves transform not only spines, but
may also transform the corresponding manifolds. Therefore, one should keep in 
memory an additional information sufficient for recovering the original
manifolds from the new ones.

Let $P$ be an almost special spine of a 
3-manifold $M$.

{\bf Move 1.} ({\it Disc replacement}). Let $P$ be a fake surface.
Suppose $D^2$ is a disc in $M$ such that  $D^2\cap P =\partial D^2$ and the
curve $\partial D^2$ is in general position in $P$. Then $D^2$ cuts of $M\setminus P$
a ball $B^3$. Let $\alpha \neq D^2$ be a 2-component of the fake surface $P\cup D^2$
such that $\alpha$ separates $B^3$ from $M\setminus B^3$. Removing from $\alpha$ an open
2-disc and collapsing the resulting polyhedron until possible, we get another 
almost special spine $P_1$ of $M$, see Fig. 2.

\begin{figure}
  \centerline{\psfig{file=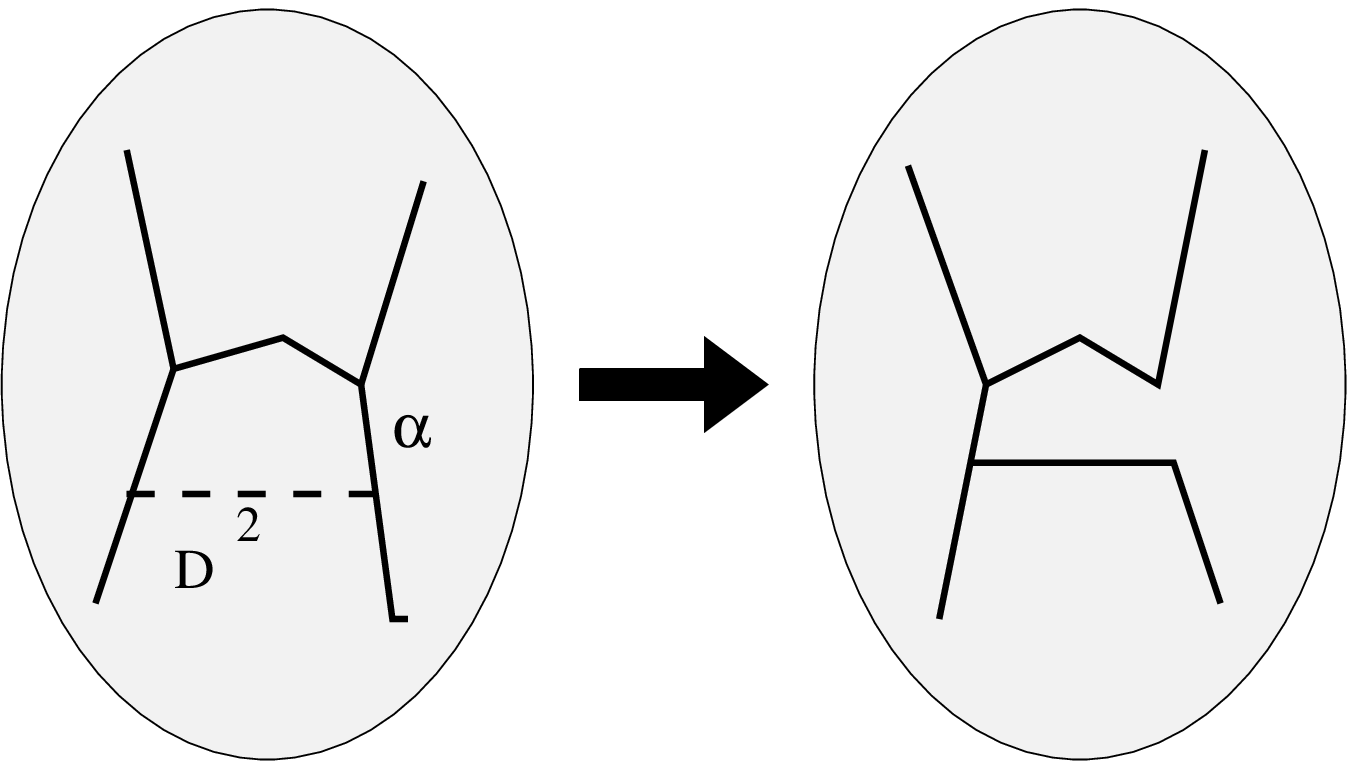,width=12cm}}  
  \caption{\label{fig2}}
  \end{figure}   

We require that Move 1 should not increase the number $v(P)$ of vertices of
$P$, i. e.
 $v(P_1)\le v(P)$.  We will say that the move is 
{\it monotone} if $v(P_1) < v(P)$, and {\it horizontal} if $v(P_1) = v(P)$.
The number $v(P\cup D^2)-v(P)$ is called {\it the 
degree} of the move.

{\bf Remark.} It is easy to see that any special spine admits only finitely
many different disc replacement moves of a given degree $n$.

In what follows we will consider disc replacement moves of degree $\le 4$.

{\bf Move 2.} ({\it Cutting a 2-component along a non-trivial circle}).
Suppose a 2-component $\alpha \subset P$ contains an orientation preserving simple
closed curve $l$ such that $l$ determines a non-trivial element of $\pi_1(M)$.
Then $\partial M$ contains at least one torus, and  there is an incompressible 
proper annulus $A\subset M$ such that $A\cap P=l$. Cutting $P$ along $l$ and 
collapsing the resulting polyhedron until possible, we get an almost special 
polyhedron $P_1\subset M$. Denote by $M_1$  a regular neighborhood  of $P_1$ in
$M$. The manifold 
$M_1$ is obtained from $M$ by cutting  along
$A$. It is easy to see that the complement   of $M_1$  in $M$ is homeomorphic
to
$Y^3=N^2\times S^1$, where $N^2$ is a disc with two holes.
Therefore, $M =  M_1\cup Y^3$, where $M_1\cap Y^3=\partial M_1\cap \partial Y^3$
consists of one  or two tori depending on whether or not the boundary circles of 
$A$ lie in different tori in $\partial M$. 

Moves 1 and 2 are basic ones. Applying them, we may obtain  almost special
spines with  one-dimensional parts as well as spines of 3-manifolds with 
several spherical
boundary components. To simplify them, we use additional moves.

{\bf Move 3.} ({\it Cutting free arcs}). Suppose $P$ contains an arc $l$ 
such that no 2-dimensional sheets are attached to $l$. Then we replace $P$ by 
$P_1=P-$Int $l$. To describe the corresponding transformation of $M$, denote by
$D^2$ a proper disc in $M$ such that 
$D^2$ intersects $l$ transversely at exactly one point.
 Then a regular neighborhood $M_1$ of $P_1$ in $M$ 
 is obtained from $M$ by cutting along $D^2$.
 In other words,
$M$ can be obtained from the new manifold $M_1$ by attaching index 
1 handle, see Fig. 3. There are three cases.

\begin{figure}
  \centerline{\psfig{file=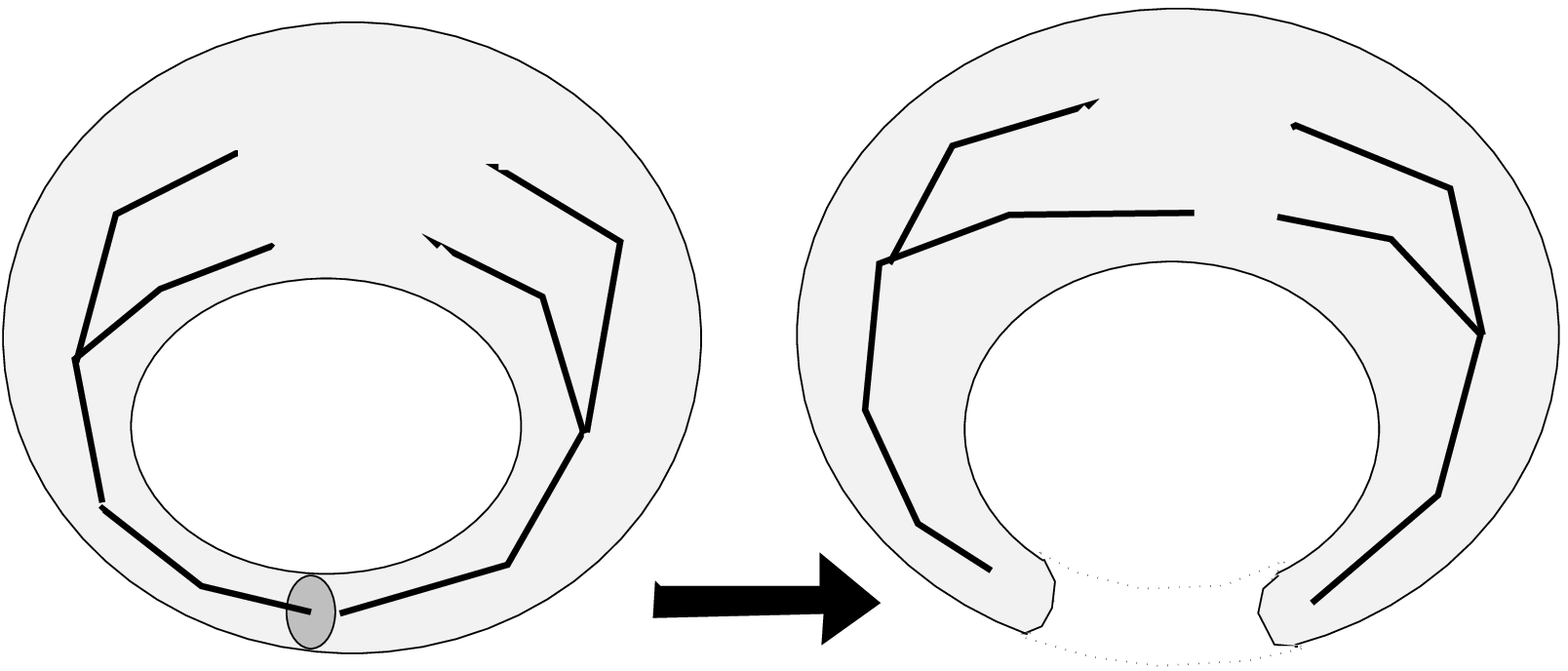,width=12cm}}  
  \caption{\label{fig3}}
  \end{figure}

{\bf A)   } $M=(M_1+D^3)   \# D^2\times S^1$, if $\partial D^2$ does not separate
$\partial M$;

{ \bf B)} $M=(M_1+D^3)  \#  S^2\times S^1$ if $D^2$ does not separate $M$, but 
$\partial D^2$ separates $\partial M$;

{ \bf C)} $M=M_1'  \# (M_1''+D^3)$ if $D^2$ separates M, where $M_1',\ M_1''$ are
the connected components of $M_1$ and $"+D^3"$ means that we fill up by a 3-ball
a spherical component of the boundary.

{\bf Move 4.} ({\it Delicate piercing}). Suppose $\partial M$ consists of at least two
components, and at least one of them is spherical. Then one can find a proper
arc $l\subset M$ such that $A=l\cap P$ is a non-singular point of $P$ 
and $l$ joints a spherical component of $\partial M$ with another one.
 Removing 
from $P$ an open disc neighborhood of $A$ and collapsing the resulting
polyhedron until possible, we get an almost special spine $P_1$ of 
$M_1=M+D^3$. We call the piercing delicate since it induces a very mild
modification of $M$.

{\bf Move 5.} ({\it Rough piercing}). Let $M$ be closed, and let $\alpha $ be a
2-component of $P$ such that the boundary curve of $\alpha$ contains
the maximal number of vertices. Removing $\alpha$ from $P$ and collapsing the
resulting polyhedron until possible, we get an almost special spine $P_1$ of
a new manifold $M_1\subset M$ such that $\partial M_1$ is a torus. Clearly,
$M\setminus $ Int$M_1$ is a solid torus, i. e. $M$ is obtained from $M_1$ by a
Dehn filling.
 
\section {Experimental results}

  We recall the notion of complexity of 3-manifolds that is naturally related
  to practically all the known methods of presenting manifolds and adequately
  describes complexity of manifolds in the informal sense of the 
  expression  [Matveev 1990].
 
 {\bf Definition.} The {\it complexity} $ c(M)$ of a compact 3-manifold  
 $M$ equals $k$ if $M$ possesses an almost special spine with $k$ vertices and
 admits no almost special spines with a smaller number of vertices.
 
    The complexity possesses the following properties:
    
    \begin{itemize}
    
    \item [ 1).]  For any integer $k$ there exists only a finite
    number of distinct closed irreducible orientable 3-manifolds of complexity
    $k$;
    
    \item [ 2).]  The complexity of the connected sum of 3-manifolds is equal
    to the sum of their complexities;
    
    \item [ 3).]  Let $M_F$ be obtained from a 3-manifold $M$ by cutting
    along a proper incompressible surface $F\subset
    M$.  Then $c(M_F)\leq c(M)$.
  
  \end{itemize}
   
 {\bf Remark.} Using moves 1, 3 and 4, one can easily prove that 
 any minimal almost special spine of a closed orientable irreducible 3-manifold
 $M$ with $c(M)> 0$ is a special one. There are exactly three 
 closed orientable irreducible   3-manifolds of complexity 0:
 $S^3, RP^3, $ and $ L(3,1)$. Their minimal almost special spines
 (a point, $RP^2$, and a fake surface without vertices) are not special.
 
Recall that a 3-manifold $M$ is a {\it graph-manifold} if it
can be obtained by pasting together several copies of  $D^2\times S^1$ 
and $Y^3=N^2\times S^1$ (where $N^2$ is a disc with two holes) along
some homeomorphisms of their boundaries.

{\bf Theorem 1.}  {\it All closed orientable 3-manifolds 
of complexity $\leq 8$ are graph-manifolds.}

 Theorem 1 was initially proved by a computer. Let us describe the main steps of
 the program.
 
 {\sc Step 1}. The computer enumerates all the special polyhedra with $\leq 8$
 vertices (there are finitely many of them);
 
 {\sc Step 2}. The computer selects spines of closed orientable 3-manifolds;
 
 {\sc Step 3}. Then it tries to apply to each spine 
  degree 4 
 disc replacement moves that
 strictly decrease the number of vertices. If such a move is possible, then
 the corresponding 3-manifold $M$ is not interesting for either it 
 has a smaller complexity
  (and we had met it earlier), or it is a connected sum of closed manifolds
  of smaller complexities. Otherwise, we go to the next step.
  
  {\sc Step 4}.  The computer  applies Move 5 (rough piercing) and 
  simplifies the new spine by Moves 1 -- 4.
   Note that  Move 5
    is allowed only if the manifold is closed or has a spherical boundary,
    and it produces a manifold with a torus boundary.

        The main observation resulting from the computer experiment is that
  if we start with a special spine of a closed orientable manifold $M$ with
  $\leq 8$ vertices, then after Moves 1 -- 5 
   we always get 
   a collection of points presenting spines of 3-dimensional spheres. 
     It means that  $M$ is a graph manifold.       
 
Later, a purely theoretical proof of Theorem 1 was found [Matveev 1990].
Note that Theorem 1 is exact in the following sense: there  exist
closed orientable 3-manifolds of complexity 9 that are hyperbolic and therefore
are not graph-manifolds. The volume of one of them is equal to 0.94...;
it is the smallest known value for volumes of closed orientable hyperbolic 
3-manifolds. [Fomenko, Matveev 1988; Hildebrand, Weeks, 1989].  

{\bf Theorem 2.} {\it If a special spine of a closed orientable
3-manifold $M$ contains  $ < 8 $ vertices and is not minimal, then it admits a degree $\leq4$
monotone disc replacement move. Any two minimal special spines of $M$ are
related by degree $\le 4$ horizontal disc replacement moves. }

   Theorem 2 had been verified by a computer program.

\section {Conjectures}
   
The following conjectures have been motivated by above-stated experimental results
as well as by other observations.

{\bf Conjecture 1.} {\it If a special spine of a compact 3-manifold is not
minimal, then the number of its vertices can be decreased by degree $\leq 4$
disc replacement moves.} 

   If the conjecture is true, then one can get a simple algorithm for
   recognition of the   unknot. Let us  apply to a spine of the   knot
   space degree $\leq 4$ disc replacement moves until possible. The knot is 
   trivial if and only if we get a circle. In the same way one can get
   a simple algorithm for recognition of the 3-sphere.
   
 Most probably, the conjecture is not true, but it is true "in general".
 In other words, the above algorithms would give   the circle or the point for
 almost all spines of the solid torus or the ball, respectively.  It means that
 we have good practical procedures for recognizing   the unknot and the
 sphere. 
 
 {\bf Conjecture 2.} {\it If a special spine of a closed graph-manifold is
 minimal, then any Move 5 (rough piercing) transforms it into a special spine of
 a graph-manifold.}
 
  The conjecture is true for all graph-manifolds up to complexity 8. It allows
  one to reduce the recognition problem for  closed graph-manifolds to 
  that for manifolds with boundaries.

   Note that if the boundary $\partial M$ of an irreducible, boundary
   irreducible graph-manifold
   $M$ is not empty, then $M$ contains an essential annulus. Conjectures 1,
   2
   and the following theorem show that Moves 1 -- 5 are "almost sufficient"
   for recognizing  graph-manifolds with boundaries and, more generally, 
   for decomposing
   3-manifolds into geometric pieces.

 {\bf Theorem 3.} {\it If a compact 3-manifold $M$ contains an essential
 annulus, then its minimal almost special spine is not special.}
 
 {\bf Remark.} If an almost special spine is not special, then it contains
 either a 2-component not homeomorphic to the 2-disc or an one-dimensional part. 
  Hence we can apply either  Move 2 or Move 3.
  
    Before proving Theorem 3, let us recall some notions of normal surface
    theory  [Haken 1961]. Let $\xi$ be a handle decomposition of a 3-manifold
    $M$ with non-empty boundary. It consists of index 0, 1, and 2  handles
called {\it balls, beams,} and {\it plates,} respectively.   
Connected components in the intersection of balls and beams are called
{\it islands}, connected components in the intersection of balls and plates are
called {\it bridges}. The boundaries of balls meet $\partial M$ along {\it
lakes}. 
    Any normal surface $F\subset M$ should intersect balls, beams and 
plates in a very specific way (see [Haken 1961]).  Particularly, 
the intersection of $F$ with balls should consist of 
{\it elementary } discs. 
The boundary curve $\partial E$ of each elementary disc $E$ 
should satisfy the following conditions:

\begin{itemize}

\item[(1)] The intersection of $\partial E$ with any 
bridge and any lake consists of no more than one segment;

\item[(2)] If $l$ is an arc in the intersection of $\partial E$ with a lake $L$
then the endpoints of $l$ should lie in different connected components of the
intersection of $L$ with islands;

\item[(3)] If a lake and a bridge are adjacent then $\partial E$ intersects 
no more than one of them.

\end{itemize}

 We will say that an elementary disc $E$  has the {\it type} $(m,n)$ 
 if the circle 
 $\partial E$ intersects $m$ bridges and $n$ lakes.
 
Any special spine $P$ of $M$ generates a handle decomposition $\xi_P$ of $M$.
Balls, beams, and  plates of the decomposition correspond to vertices,
edges and 2-components of $P$, respectively. The boundary of each ball contains exactly
four islands, and any two of them are joined by exactly one bridge, see Fig. 4.

\begin{figure}
  \centerline{\psfig{file=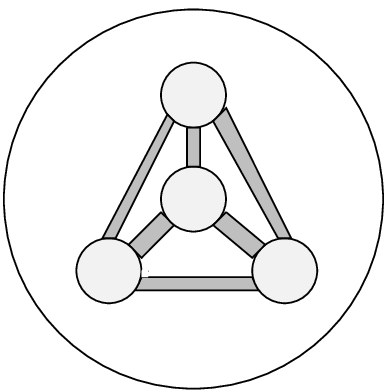,width=6cm}}  
  \caption{\label{fig4}}
  \end{figure}   

It is not hard to see that any elementary disc  for $\xi(P)$ has one of 
the following types: $(4,0), (3,0), (2,1), (1,2), (0,2), (0,3), (0,4) $. 
 Each type determines 
the corresponding elementary disc in a unique way (up to homeomorphisms 
of the ball
taking islands to islands, bridges to bridges, and lakes to lakes), except 
the  type (0,3)  that determines two elementary discs.

For any beam $D^2\times I$ (with $D^2\times \{0\}$ and $D^2\times \{1\}$
being islands), the disc $D^2\times \{1/2\}$ is called the {\it transverse
disc} of the beam.

{\bf Definition.} Let $A$ be a proper annulus in a 3-manifold $M$ 
with a special spine $P$ such that $A$ is normal with respect to $\xi_P$.
We say that $A$ has {\it a tail} if the intersection of $A$ with the 
transverse disc of a beam
contains a proper arc $l$ such that the endpoints of $l$ lie in the same
circle of $\partial A$. The arc $l$ cuts off $A$ a disc $D_l$. We will
refer to $D_l$ as to a {\it tail} of $A$.

   {\bf Lemma 1.}  {\it If the generated by a special spine $P$ of a 3-manifold
   $M$ handle decomposition $\xi_P$ contains a normal annulus $A$ with a tail
 $D_l$, then $P$ is not minimal.}

    {\bf Proof.}  Denote by $M_{D_l}$ the  3-manifold obtained from $M$ by
    cutting along $D_l$. Evidently, $M_{D_l}$ is homeomorphic with $M$. 
    The tail decomposes the balls of $\xi_P$ into balls,
    plates into plates, and beams into beams except the beam $B_0$ containing $l$.
Coherently collapsing new balls, beams, and plates onto 2-dimensional subsets,
we get an almost special spine $P'$ of $M_{D_l}$ . Since each ball of $\xi_P$
contains no more than one vertex of $P'$, we have $v(P')\leq v(P)$, where
$v(P)$ denotes the number of vertices. Note  that $P'$ has a free edge arising
from
cutting and collapsing the beam $B_0$, see Fig. 5. After collapsing $P'$
through this free edge, we get an almost special spine of  $M_{D_l}$ with
a fewer number of vertices.

\begin{figure}
  \centerline{\psfig{file=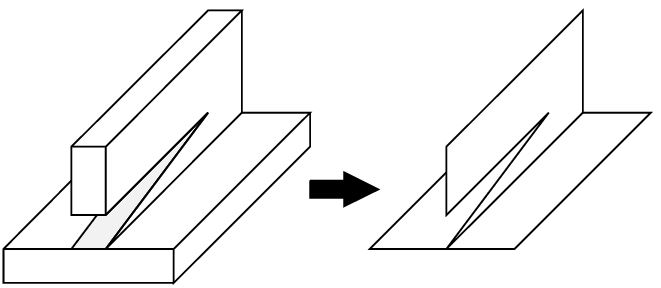,width=12cm}}  
  \caption{\label{fig5}}
  \end{figure}   
  
{\bf Proof of Theorem 3.} Let $P$ be a special spine of a 3-manifold $M$
with an essential annulus. Replace the annulus by an annulus $A$ that
is normal with respect to the generated by $P$ handle decomposition $\xi_P$ of $M$.
If $A$ has a tail, then we  apply Lemma 1 to find a simpler spine of 
$M$. Assume that $A$ has no tails. Since each elementary disc of type
$(0,3)$ or $(0,4)$ in $A$ would determine at least one tail,   only types
 $(4,0), (3,0), (2,1), (1,2), (0,2)$ for elementary discs in $A$ are possible.    
Moreover, if $E$ is an elementary disc of type $(1,2)$ or $(0,2)$, then
two arcs in $\partial E \cap \partial M$ must lie in different components of
$\partial A$.

   Let us cut now $M$ along $A$ such that one component $S$ of $\partial A$ is
   preserved. In other words, we remove from $M$ the subset
   $S^1\times (0,1]\times I$, where $A=S^1\times [0,1]$ and 
   $S^1\times [0,1]\times I $ is a thin regular neighborhood of $A$ in $M$.
As above,  coherently collapsing the new balls, beams, and plates onto 
2-dimensional subsets, we get an almost special spine $P'$ of $M$ with
$v(P')\leq v(P)$. Moreover, if at least one elementary
disc of  type $(1,2)$ is present, then  $v(P')< v(P)$. It is because
each type $(1,2)$ elementary disc in the intersection of $A$ with a ball of
$\xi_P$ annihilates the corresponding vertex of $P$, see Fig. 6.  
  \begin{figure}
  \centerline{\psfig{file=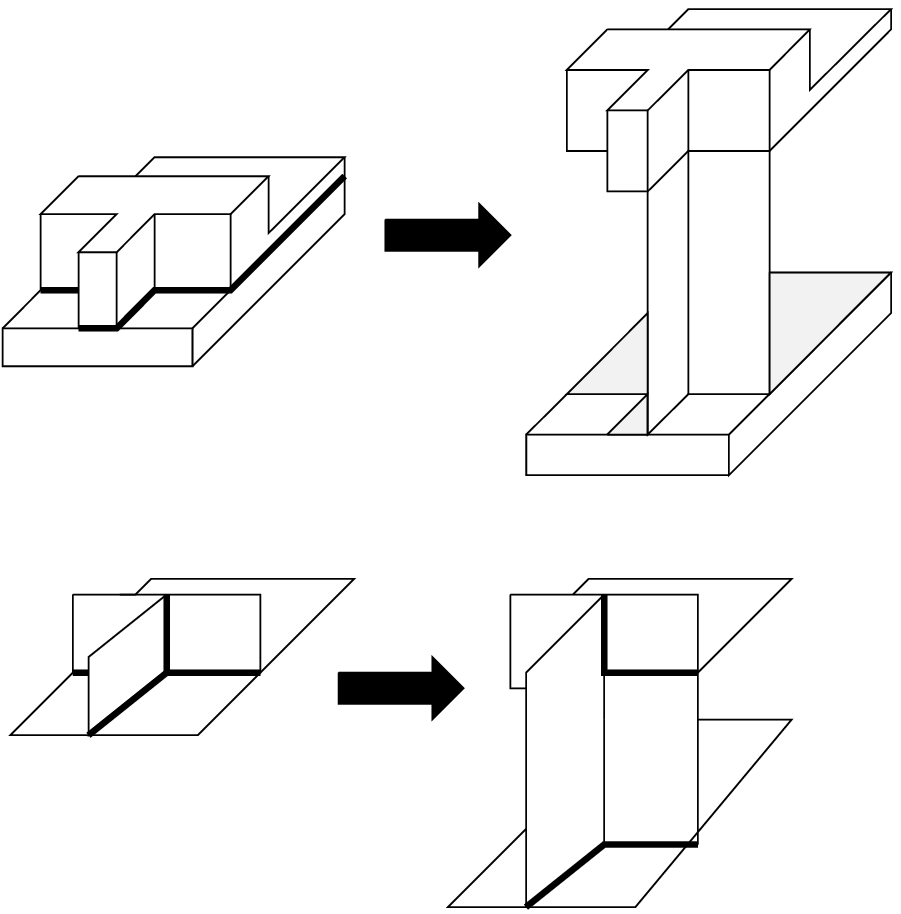,width=12cm}}  
  \caption{\label{fig6}}
  \end{figure}   

   We conclude the proof with the following remark: if there are no type 
   $(1,2)$ elementary discs in $A$, then all elementary discs in $A$ 
   have type $(0,2)$. In this case the core circle of $ A$ can be 
   shifted into a  2-component of $P$.  Since all 2-components of $P$
   are 2-cells, it implies that the core circle is contractible, which
   contradicts  the assumption that $A$ is incompressible.

\section{The algorithm}

  Let $M$ be a compact 3-manifold such that $\partial M$ is either empty or
  consists  of  tori. Our goal is to decompose $M$ into 
  geometric pieces,
  particularly, to determine whether or not $M$ is a graph-manifold.
  
  {\sc Step 1}. Construct a special spine $P$ of $M$;
  
  {\sc Step 2}. Apply to $P$  Moves 1 -- 4 until possible. In the case of Move 2
  when $Y^3\times S^1$ is cut off we store the information how  
  $Y^3\times S^1$ is attached to the remaining part of $M$. It can be done by
  selecting a meridian-longitude pair on boundary tori and controlling their
  behavior under further moves. 
  
  {\sc Step 3}. Assume that we have started with a non-closed manifold  $M$.
 Then by Theorem 3 (and  modulo Conjecture 1), $M$ is a graph-manifold if and only if after
 Step 2 we get a collection of points. We use the stored information to
 present  $M$ as a union of Seifert manifolds with explicitly given
 parameters and gluing matrices.
 
 {\sc Step 4}. In general, we get a collection of special spines of irreducible,
 boundary irreducible manifolds. Then we apply Move 5 to those of them 
 that are spines of closed manifolds. 
  
{\sc Step 5}.  Iterate Steps 2 -- 4 until possible. 

     If we have started with a graph-manifold  $M$, then (modulo 
     Conjectures 1 and 2) we get an explicit presentation of $M$
     as a connected sum of unions of Seifert manifolds with known
     parameters and known gluing matrices. In general, we may get
     unknown pieces that should be tested for hyperbolicity by comparing 
     with J. Weeks table. 

\small

\section*{REFERENCES}

\begin{description}

\item [{[Fomenko, Matveev 1988]}]
 S. V. Matveev,  A. T. Fomenko,    
``Isoenergetic  surfaces   of   
    Hamiltonian systems,  enumeration  of three-dimensional manifolds in 
   increasing order of
     complexity,  and    computation  of volumes of   closed
    hyperbolic manifolds'' (Russian),  {\it Uspekhi  Mat.  Nauk } 
    {\bf 43}  (1988), no. 1(259), 5-22; English transl. in 
    {\it Russian Math. Surveys} {\bf 43} (1988), no. 1, 3--24.

\item [{[Haken 1961]}] W.  Haken, ``Theorie  der Normalfl\"{a}chen. Ein Isotopiekriterium
   f\"{u}r  der Kreisknoten'', {\it Acta Math.,} {\bf 105} (1961),  245--375.

\item [{[Hildebrand, Weeks 1989]}] M. Hildebrand, J. Weeks, ``A computer generated
census of cusped hyperbolic 3-manifolds'', {\it Computers and mathematics},
(Cambridge, MA, 1989), 53--59, Springer, New York-Berlin, 1989.

\item [{[Matveev 1990]}] S. V. Matveev,
``Complexity theory of 3-manifolds'',  {\it Acta Applicandae  Mathematicae,}
{\bf 19} (1990), 101--130.
    
\item [{[Scott 1983]}] P. Scott, ``The geometries of 3-manifolds'', {\it Bull.
London Math. Soc.,} {\bf 15} (1983), 401--487.

\item [{[Waldhausen 1967]}] F. Waldhausen, ``Eine Klasse von 3-dimensionalen
Mannigfaltigkeiten. I, II'', {\it Invent. Math.} {\bf 3} (1967), 308--333; ibid.
{\bf 4} (1967), 87--117.

\end{description}

\end{document}